\documentclass[review]{amsart}
\usepackage{amsmath,amsfonts,amsthm,amsopn,color,amssymb,enumitem}
\usepackage{palatino}
\usepackage{graphicx}
\usepackage[english]{babel}
\usepackage{caption}
\usepackage{subcaption}
\usepackage[colorlinks=true]{hyperref}
\hypersetup{urlcolor=blue, citecolor=red, linkcolor=blue}
\usepackage[utf8]{inputenc}
\usepackage{mathscinet}

\usepackage{cleveref}
\crefname{equation}{}{}
\crefname{figure}{}{}

\usepackage{esint}
\usepackage{color}

\usepackage[title]{appendix}

\usepackage[backend=bibtex,style=numeric,sorting=nyt,firstinits=true,isbn=false,url=false,doi=false]{biblatex}
\usepackage{bm}
\numberwithin{equation}{section}

\newtheorem{thm}{Theorem}[section]
\newtheorem{defn}[thm]{Definition}
\newtheorem{lemma}[thm]{Lemma}
\newtheorem{cor}[thm]{Corollary}
\newtheorem{re}[thm]{Remark}

  \newcommand{\dist}{\mbox{dist}}
  \newcommand{\diam}{{\mbox{diam}}}

  \numberwithin{equation}{section}
  \numberwithin{figure}{section}

\begin{document}

\title[Boundary Regularity on ${C^{1,\mathrm{Dini}}}$ Domains]{Boundary Regularity for Fully Nonlinear Parabolic equations on $\bm{C}^{1,\mathrm{Dini}}$ Domains}

\author{Jiqi Dong}
\address{J. Dong: School of Mathematical Sciences, University of Jinan, Jinan, P.R.China 250022}
\email{17860357019@163.com}

\author{Xuemei Li}
\address{X. Li: School of Mathematical Sciences, University of Jinan, Jinan, P.R.China 250022}
\email{sms\_lixm@ujn.edu.cn}

\author{Yuanyuan Lian}
\address{Y. Lian: Departamento de An\'{a}lisis Matem\'{a}tico, Universidad de Granada, Campus Fuentenueva, 18071 Granada, Espa\~{n}a}
\email{lianyuanyuan.hthk@gmail.com; yuanyuanlian@correo.ugr.es}

\date{\today}


\thanks{X. L. has been supported by National Natural Science Foundation of China (No.12401257) and the Natural Science Foundation of Shandong Province, China (No. ZR2024QA045). Y. L. has been supported by the Grants PID2020-117868GB-I00 and PID2023-150727NB-I00 of the MICIN/AEI}

\begin{abstract}
We establish the boundary pointwise Lipschitz regularity on exterior $C^{1,\mathrm{Dini}}$ domains and the Hopf lemma on interior $C^{1,\mathrm{Dini}}$ domains for fully nonlinear parabolic equations by a unified perturbation method. In fact, above two regularity hold for more general solution sets, i.e., the Pucci's class $S^*(\lambda, \Lambda, f)$. Furthermore, based on the boundary pointwise Lipschitz regularity, we obtain the global ${W}^{2,\delta}$ regularity on exterior $C^{1,\mathrm{Dini}}$ domains for any $0<\delta<1$, which is new even for the harmonic functions.
\end{abstract}

\maketitle

\section{introduction}\label{Intro}
We study viscosity solutions of fully nonlinear uniformly parabolic equations
\begin{equation}\label{e.Dirichlet}
\begin{cases}
u_t-F(D^2u)=f& ~~\mbox{in}~~\Omega\cap Q_1,\\
u=g& ~~\mbox{on}~~\partial \Omega\cap Q_1,
\end{cases}
\end{equation}
where $\Omega \subset \mathbb{R}^{n+1}$ is a bounded domain and $Q_1$ is the unit parabolic cube (see \Cref{pre} for the notations and preliminaries). In this paper, we establish the boundary regularity on $C^{1,\mathrm{Dini}}$ domains.

Boundary regularity for fully nonlinear parabolic equations has attracted much attention during last decades. In 1983, Krylov \cite{MR688919} gave the boundary $C^{1,\alpha}$ and $C^{2,\alpha}$ a priori estimates. Later, Wang \cite{MR1139064} proved the boundary pointwise $C^{1,\alpha}$ regularity. The third author and Zhang \cite{LZ_2022} established boundary pointwise $C^{k,\alpha}$ regularity for any $k\geq 1$. Adimurthi, Banerjee and Verma \cite{MR4073515} obtained the boundary $C^{2,\alpha}$ regularity (not pointwise) by flattening the boundary.

All regularity above require that the domains are $C^{1,\alpha}$ or smoother. On $C^{1,\mathrm{Dini}}$ (weaker than $C^{1,\alpha}$) domains, by constructing proper barriers, Kamynin and Khimchenko proved the boundary Lipschitz regularity \cite{MR437947,MR427828,KK2} and the Hopf lemma \cite{MR333453,KK} respectively. However, they only proved the results for smooth solutions of linear equations. Moreover, some additional assumptions on the Dini modulus were proposed. We also refer to \cite{Safonov2008,MR3167627,MR4713521,torreslatorre2025} for the corresponding results for elliptic equations on $C^{1,\mathrm{Dini}}$ domains.

Up to our knowledge, there is no boundary Lipschitz regularity and Hopf lemma for viscosity solutions of fully nonlinear parabolic equations on $C^{1,\mathrm{Dini}}$ domains. In this paper, following the idea of \cite{MR4713521}, we fill this gap.

For Cadel\'{o}n-Zygmund estimates, Wang \cite{MR1135923} established interior $W^{2,\delta}$ estimates for Pucci's class $S(\lambda,\Lambda,f)$ and then established interior and global $W^{2,p}$ estimates for \eqref{e.Dirichlet} on $C^{1,1}$ domains. Recently, the second author \cite{LI2025109459} established global $W^{2,\delta}$ regularity (for some small $\delta$) for $S(\lambda,\Lambda,f)$ on $C^{1,\alpha}$ domains. In this paper, with the aid of the boundary pointwise Lipschitz regularity, we prove the global $W^{2,\delta}$ regularity for any $0<\delta<1$ on exterior $C^{1,\mathrm{Dini}}$ domains. Note that this regularity is new even for the harmonic functions.

It is worth pointing out that the pointwise regularity is more essential than the classical local and global regularity since we can obtain the later from the former but not vice versa. In addition, the pointwise regularity has some distinctive applications that can not be fulfilled by the classical local regularity. The following are several examples. Figalli and Shahgholian \cite[Proof of Theorem 1.2]{MR3198649} utilized the interior pointwise $C^{2,\alpha}$ regularity to prove the optimal $C^{1,1}$ regularity for solutions of obstacle type free boundary problems. The boundary pointwise $C^{k,\alpha}$ ($k\geq 1$) regularity can be used to give a direct and simple proof of the higher regularity of free boundaries (see \cite{MR4644419} and \cite{LZ_2022}) without using the partial hodograph-Legendre transformation as in \cite{MR0440187}. Li, the second author and Zhang \cite{MR4631039} employed the boundary pointwise $C^{1,\alpha}$ regularity to establish the global $W^{2,p}$ regularity   on $C^{1,\alpha}$ domains. The third author \cite{lian2024nodalsets} used the interior pointwise regularity to prove the higher order smoothness of the nodal set of solutions of elliptic and parabolic equations.

\medskip

Now, our main results are stated as follows. For boundary pointwise regularity, we always assume that $0:=(0,0)\in\partial \Omega$ and study the pointwise regularity at $0$. The first result concerns the \emph{boundary pointwise Lipschitz regularity} on exterior $C^{1,\mathrm{Dini}}$ domains (we refer to \Cref{pre} for the notions and notations).
\begin{thm}[\textbf{Boundary Lipschitz regularity}]\label{t3.1}
Let $u$ be a viscosity solution of
	\begin{equation}\label{3.1}
		\begin{cases}
          u\in S^*(\lambda,\Lambda,f) &~~\mathrm{in}~~\Omega\cap Q_1,\\
		  u=g &~~\mathrm{on}~~\partial\Omega\cap Q_1.
		\end{cases}
	\end{equation}
Suppose that $\Omega$ satisfies the exterior $C^{1,\mathrm{Dini}}$ condition at $0\in \partial \Omega$ and
\begin{equation*}
  f\in C^{-1,\mathrm{Dini}}(0),\quad g\in C^{1,\mathrm{Dini}}(0).
\end{equation*}
Then $u\in C^{0,1}(0)$, i.e.,
	\begin{equation}\label{3.2}
		|u(x,t)-u(0)|\leq C|(x,t)|
		(\|u\|_{L^{\infty}(\Omega\cap Q_1)}+\|f\|_{C^{-1,\mathrm{Dini}}(0)}+\|g\|_{C^{1,\mathrm{Dini}}(0)}),\
		\forall (x,t)\in \ \Omega\cap Q_1,
	\end{equation}
where $C$ depends only on $n$, $\lambda$, $\Lambda$, $\omega_\Omega$, $\omega_f$ and $\omega_g$.
Here, $\omega_{\Omega}$ etc. denote the Dini modulus corresponding to $\Omega$ etc. (see \Cref{pre} for the precise definitions).
\end{thm}	

\begin{re}\label{re1.1}
Note that the Pucci's class is more general than an equation like \eqref{e.Dirichlet}. There exists a function $u \in S(\lambda,\Lambda,0)$ which cannot be a viscosity solution of $u_t-F(D^2u)=0$ (see \cite[Section 1.7.2]{MR4420058}).
\end{re}

\begin{re}\label{re.b.Lip}
To obtain the boundary pointwise Lipschitz regularity, the exterior $C^{1,\mathrm{Dini}}$ condition on the domain is optimal (see \cite[Theorem 1.15]{MR4713521} for the elliptic  equation which is a special case of parabolic equations).
\end{re}

\begin{re}\label{re.g.3}
For the Pucci's class $S(\lambda, \Lambda, f)$, the optimal interior regularity is the $C^{\alpha}$ regularity (see \cite[Theorem 4.19]{MR1135923}). However, we can obtain higher boundary regularity.
\end{re}

\begin{re}\label{re.b.Lip2}
This boundary pointwise Lipschitz regularity also holds for \textbf{parabolic normalized $p$-Laplace equations} because it belongs to $S^*(\lambda,\Lambda,f)$ (see \cite[Theorem 1.3 and Theorem 1.4]{MR4656656}). Similarly, the Hopf lemma (\Cref{t3.1'}) and global $W^{2,\delta}$ regularity (\Cref{t4.1}) also hold for parabolic normalized $p$-Laplace equations.
\end{re}

\begin{re}\label{re1.6}
It is standard (see \cite[Section 2]{MR1135923}) that for a parabolic domain $\Omega$, its parabolic boundary $\partial \Omega$ can be decomposed as $\partial \Omega=\partial_x \Omega\cup \partial_b \Omega$, where $\partial_x \Omega$ denotes the lateral boundary and $\partial_b \Omega$ denotes the bottom boundary. From the definition of $C^{1,\mathrm{Dini}}$ condition on the domain (see \Cref{de1.2} and \Cref{de1.3}), we know that $0$ must belong to the lateral boundary $\partial_x\Omega$. Hence, \Cref{t3.1} is a regularity on the lateral boundary.
\end{re}

By combining with the interior  $C^{1}$ regularity of \eqref{e.Dirichlet}, we have the following \emph{global Lipschitz regularity}.

\begin{cor}\label{cor1.1}
Let $u$ be a viscosity solution of \eqref{e.Dirichlet} and $\Omega$ be an exterior $C^{1,\mathrm{Dini}}$ domain. Suppose that
\begin{equation*}
  f\in C^{-1,\mathrm{Dini}}(\overline\Omega), \quad
  g\in C^{1,\mathrm{Dini}}(\partial\Omega).
\end{equation*}
Then for any $\Omega'\subset \Omega$ with $\mathrm{dist}(\partial_b \Omega,\Omega')>0$, we have
 $u\in C^{0,1}(\overline {\Omega'})$ and
\begin{equation*}
  \|u\|_{C^{0,1}(\overline {\Omega'})}\leq C\left(\|u\|_{L^{\infty}(\Omega)}+\|f\|_{C^{-1,\mathrm{Dini}}(\overline\Omega)}
  +\|g\|_{C^{1,\mathrm{Dini}}(\partial\Omega)}\right),
\end{equation*}
where $C$ depends only on $n,\lambda,\Lambda$, $\omega_{\Omega}$, $\omega_f$, $\omega_g$ and $\mathrm{dist}(\partial_b \Omega,\Omega')$. Here, $\mathrm{dist}(\partial_b \Omega,\Omega')$ denotes the distance between the bottom boundary $\partial_b \Omega$ and $\Omega'$.
\end{cor}

\begin{re}\label{re.g.1}
With the condition $f\in C^{-1,\mathrm{Dini}}$, one can get the interior $C^{1}$ regularity for \eqref{e.Dirichlet} by the perturbation argument. However, we can not find a direct reference to this result. We refer to \cite{lian2020pointwise} for an elliptic counterpart.
\end{re}

In addition, with the aid of the higher interior regularity for $S^*(\lambda,\Lambda,f)$ if $\lambda$ and $\Lambda$ are enough close, we have the following corollary.

\begin{cor}\label{cor1.2}
Let $u$ be a viscosity solution of \eqref{3.1} and $\Omega$ be an exterior $C^{1,\mathrm{Dini}}$ domain. Suppose that
\begin{equation*}
  f\in C^{-1,\mathrm{Dini}}(\overline\Omega), \quad
  g\in C^{1,\mathrm{Dini}}(\partial\Omega).
\end{equation*}
Assume $\Lambda\leq (1+\eta) \lambda$, where $0<\eta<1$ is a small constant depending only on the dimension $n$.

Then for any $\Omega'\subset \Omega$ with $\mathrm{dist}(\partial_b \Omega,\Omega')>0$, we have
 $u\in C^{0,1}(\overline {\Omega'})$ and
\begin{equation*}
  \|u\|_{C^{0,1}(\overline {\Omega'})}\leq C\left(\|u\|_{L^{\infty}(\Omega)}+\|f\|_{C^{-1,\mathrm{Dini}}(\overline\Omega)}
  +\|g\|_{C^{1,\mathrm{Dini}}(\partial\Omega)}\right),
\end{equation*}
where $C$ depends only on $n,\lambda,\Lambda$, $\omega_{\Omega}$, $\omega_f$, $\omega_g$ and $\mathrm{dist}(\partial_b \Omega,\Omega')$.
\end{cor}

\begin{re}\label{re.g.2}
As mentioned in \Cref{re.g.1}, one can obtain the interior $C^{0,1}$ regularity for \eqref{3.1} if $\lambda$ and $\Lambda$ are quite close.
\end{re}

\begin{re}\label{re.b.Lip3}
By H$\ddot{o}$lder inequality, the condition $f\in L^p(\Omega)$ for  $p>n+2$ implies $f\in C^{-1,\mathrm{Dini}}(\overline\Omega)$. Thus, \Cref{cor1.1} and \Cref{cor1.2} hold with  $f\in L^p(\Omega)$ for $p>n+2$.
\end{re}

Next, we show the \emph{Hopf lemma} on interior $C^{1,\mathrm{Dini}}$ domains.
\begin{thm}[\textbf{Hopf lemma}]\label{t3.1'}
Let $u\in S^*(\lambda,\Lambda,f)$ in $\Omega\cap Q_1$ with $u\geq 0$ and $u(0)=0$. Suppose that $\Omega$ satisfies the interior $C^{1,\mathrm{Dini}}$ condition at $0\in \partial \Omega$ and $f\in C^{-1,\mathrm{Dini}}(0)$ with
\begin{equation*}
  \|f\|_{C^{-1,\mathrm{Dini}}(0)}\leq \varepsilon_0\cdot u(e_n/2,-1/2),
\end{equation*}
where $0<\varepsilon_0<1$ is a small constant depending only on $n,\lambda,\Lambda$ and $\omega_\Omega$.

Then for any $l=(l_1,l_2,...,l_n)\in \mathbb{R}^n$
with $|l|=1$ and $l_n> 0$,
\begin{equation}\label{hf}
	u(rl,t)\geq c u(e_n/2,-1/2) l_nr,\quad \forall~0<r<\eta, -\eta^2<t\leq 0.
\end{equation}
where $c$ depends only on $n,\lambda,\Lambda$, $\omega_\Omega$ and $\eta$ depends also on $l$.
\end{thm}

\begin{re}\label{re.b.Hop}
For the assumption of Hopf Lemma on domains, the interior $C^{1,\mathrm{Dini}}$ condition is the optimal (see \cite[Theorem 1.18]{MR4713521}). Up to our knowledge, above theorem is the first Hopf lemma concerning non-homogenous righthand term $f$. Hence, above theorem is new even for the Poisson's equation.
\end{re}

\begin{re}\label{re1.2}
Note that the smallness condition on $f$ can not be removed. A typical counterexample is the following:
$u(x)=x_n^2$ is a solution of
	\begin{equation*}
		\begin{cases}
         \Delta u=2&~~\mathrm{in}~~B_1^+,\\
		  u=0 &~~\mathrm{on}~~T_1:=\left\{(x',x_n): |x'|<0, ~~ x_n=0\right\}.
		\end{cases}
	\end{equation*}
Clearly, the Hopf lemma fails for $u$.
\end{re}

%

The last main result concerns the \emph{global $W^{2,\delta}$ regularity} on exterior $C^{1,\mathrm{Dini}}$ domains.

\begin{thm}[\textbf{Global $W^{2,\delta}$ regularity}]\label{t4.1}
Let $F$ be convex or concave in $M$ and $u$ be a viscosity solution of \eqref{e.Dirichlet}. Suppose that $\Omega$ is both a Lipschitz domain and an exterior $C^{1,\mathrm{Dini}}$ domain. Assume that
\begin{equation*}
  f\in  L^p(\Omega) \ (p>n+2)  \quad \mathrm{and}\quad  g\in C^{1,\mathrm{Dini}}(\partial \Omega).
\end{equation*}
Then for any $0<\delta<1$ and $\Omega'\subset \Omega$ with $\mathrm{dist}(\partial_b \Omega,\Omega')>0$, we have $u\in W^{2,\delta}(\Omega')$ and
\begin{equation}\label{1.2}
	\|u \|_{ W^{2,\delta}(\Omega')}\leq C\left( \|u\|_{L^{\infty}(\Omega)}
   +\|f\|_{L^{p}(\Omega)}
   +\|g\|_{C^{1,\mathrm{Dini}}(\partial \Omega)}\right),
\end{equation}
where $C$ depends only on $n$, $\lambda$, $\Lambda$, $p$, $\delta$, $\omega_{\Omega},\omega_g$ and $\mathrm{dist}(\partial_b \Omega,\Omega')$.
\end{thm}

\begin{re}\label{re.w2d2}
Here we need $p>n+2$ since we need the boundary Lipschitz estimate during the proof (see \Cref{re.b.Lip3}).
\end{re}


\begin{re}\label{re.w2d0}
For the Pucci's class $S^*(\lambda, \Lambda, f)$ (problem \eqref{3.1}), we only can obtain the global $W^{2,\delta}$ regularity where $\delta$ is very small (e.g. \cite{LI2025109459}). The reason is that for $S^*(\lambda, \Lambda, f)$, we only have the interior $W^{2,\delta_0}$ regularity for some small $\delta_0$ (see \cite[Theorem 4.11]{MR1135923}).
\end{re}

As far as we know, \Cref{t4.1} is the first $W^{2,\delta}$ regularity with $\delta$ close to $1$ arbitrarily (rather than quite small) on exterior $C^{1,\mathrm{Dini}}$ domains. Hence, this result is new even for the heat equation and we present it as follows.
\begin{cor}\label{Lap_W2d}
Let $u$ be a solution of
\begin{equation*}
\begin{cases}
u_t-\Delta u=f& ~~\mbox{in}~~\Omega,\\
u=g& ~~\mbox{on}~~\partial \Omega.
\end{cases}
\end{equation*}
Suppose that $\Omega$ is both a Lipschitz domain and an exterior $C^{1,\mathrm{Dini}}$ domain,
\begin{equation*}
  f\in  L^p(\Omega) \ (p>n+2)  \quad \mathrm{and}\quad  g\in C^{1,\mathrm{Dini}}(\partial \Omega).
\end{equation*}
Then for any $0<\delta<1$ and $\Omega'\subset \Omega$ with $\mathrm{dist}(\partial_b \Omega,\Omega')>0$, we have $u\in W^{2,\delta}(\Omega')$ and
\begin{equation*}
	\|u \|_{ W^{2,\delta}(\Omega')}\leq C\left( \|u\|_{L^{\infty}(\Omega)}
   +\|f\|_{L^{p}(\Omega)}
   +\|g\|_{C^{1,\mathrm{Dini}}(\partial \Omega)}\right),
\end{equation*}
where $C$ depends only on $n$, $\lambda$, $\Lambda$, $p$, $\delta$, $\omega_{\Omega}$, $\omega_g$ and $\mathrm{dist}(\partial_b \Omega,\Omega')$.
\end{cor}

\begin{re}\label{re.Lap1}
Even for the elliptic case, the global $W^{2,\delta}$ regularity for any $0<\delta<1$ on exterior $C^{1,\mathrm{Dini}}$ domains is new.
\end{re}

We illustrate our idea as follows. For the boundary pointwise Lipschitz regularity, motivated by \cite{MR4713521}, we regard the curved boundary as the perturbation of a flat boundary. The proof of the Hopf lemma also employs the same perturbation idea. For the global $W^{2,\delta}$ estimate, traditional method of flattening the boundary is invalid for the $C^{1,\mathrm{Dini}}$ domains. Inspired by \cite{MR4583741} and \cite{LI2025109459},
with the aid of the boundary pointwise $C^{0,1}$ estimate and the Whitney decomposition, we derive the global $W^{2,\delta}$ estimate for any $0<\delta<1$ from the known interior $W^{2,p}$ estimate (see \cite[Theorem 5.6]{MR1135923})

The paper is organized as follows. In \Cref{pre}, some preliminaries, in particular the pointwise characterizations of functions and domains, are introduced. In \Cref{Lip_Hopf}, we prove \Cref{t3.1} and \Cref{t3.1'} for $S^*(\lambda,\Lambda,f)$ on the exterior and the interior $C^{1,\mathrm{Dini}}$ domains respectively. Based on the boundary pointwise $C^{0,1}$ estimate, we prove \Cref{t4.1} via the Whitney decomposition in \Cref{W2d}.

\section{Preliminaries}\label{pre}

In this section, we introduce some notations and notions, in particular, the pointwise characterizations of functions and domains. Let us first recall the parabolic terminology. Throughout this paper, we denote $\{e_i\}^{n}_{i=1}$ as the standard basis of $\mathbb{R}^{n}$, i.e., $e_i=(0,...0,\underset{i^{th}}{1},0,...0)$. In addition, denote $x'=(x_1,x_2,...,x_{n-1})\in \mathbb{R}^{n-1}$,
$x=(x',x_n)\in \mathbb{R}^{n}$ and $(x,t)=(x',x_n,t)\in \mathbb{R}^{n+1}$. The parabolic norm of $(x,t)$ is defined as
	\begin{equation*}
|(x,t)|=	\left\{
	\begin{aligned} &(|x|^2+|t|)^{1/2}\ \    &t\leq 0,\\
		 &\infty\ \   &t>0.
	\end{aligned}
	\right.
\end{equation*}
Moreover, let
\begin{equation*}
  B_r(x_0)=\{x\in \mathbb{R}^{n}:|x-x_0|<r\}~~\mbox{ and }~~Q_r(x_0,t_0)=B_r(x_0)\times(t_0-r^2,t_0].
\end{equation*}
Denote $B_r(0)$ and $Q_r(0)$ as $B_r$ and $Q_r$ respectively for simplicity. The derivatives of a function is defined as
\begin{equation*}
  D^k f=\{D^\sigma_xD^\gamma_t f: |\sigma|+2\gamma=k\},\quad k\geq 0
\end{equation*}
and its norm is defined as
\begin{equation*}
  |D^k f|=\left(\sum_{|\sigma|+2\gamma=k}|D^\sigma_xD^\gamma_t f|^2\right)^{1/2},
\end{equation*}
where $\gamma\in \mathbb{N}$ and $\sigma\in \mathbb{N}^n$ is the standard multi-index.

\medskip

In this paper, we consider viscosity solutions, which are defined as follows and we refer to \cite{MR1351007, MR1789919} and \cite{MR1135923} for more details.
\begin{defn}\label{d-viscoF}
Let $p\geq n+1$ and $f\in L^p(\Omega)$. We say that $u\in C(\Omega)$ is an $L^p$-viscosity subsolution (resp., supersolution) of \eqref{e.Dirichlet} if
\begin{equation*}
  \begin{aligned}
   & ess~\underset{(y,s) \to (x,t)}{\lim\inf}\left(\varphi_t (y,s)-F(D^2\varphi (y,s))-f(y,s) \right)\leq 0\\
   & \left(\mathrm{resp.},~ess~\underset{(y,s) \to (x,t)}{\lim\sup}\left(\varphi_t (y,s)-F(D^2\varphi (y,s))-f(y,s) \right)\geq 0\right)
  \end{aligned}
\end{equation*}
provided that for $\varphi\in W^{2,p}(\Omega)$, $u-\varphi$ attains the local maximum (resp., minimum) at $(x,t)\in\Omega$.

We call $u\in C(\Omega)$ an $L^p$-viscosity solution of \eqref{e.Dirichlet} if it is both an $L^p$-viscosity subsolution and supersolution of \eqref{e.Dirichlet}.
\end{defn}

\begin{re}\label{re2.1}
In this paper, we only consider the $L^{n+1}$-viscosity solution and a solution always means the $L^{n+1}$-viscosity solution.
\end{re}

Throughout this paper, we always assume that the fully nonlinear operator $F$ is uniformly elliptic with positive constants $\lambda$ and $\Lambda$, that is, for any $M,N\in \mathcal{S}^n$,
\begin{equation}\label{e.Stru}
\mathcal{M}^-(M,\lambda,\Lambda)\leq F(M+N)-F(N)
\leq \mathcal{M}^+(M,\lambda,\Lambda),
\end{equation}
where $\mathcal{S}^n$ denotes the set of $n\times n$ symmetric matrices; the $\mathcal{M}^-$ and $\mathcal{M}^+$ are the Pucci's extremal operators, i.e.,
\begin{equation*}
\mathcal{M}^-(M,\lambda,\Lambda):=\lambda\sum_{\lambda_i>0} \lambda_i+\Lambda\sum_{\lambda_i<0} \lambda_i,\quad
\mathcal{M}^+(M,\lambda,\Lambda):=\Lambda\sum_{\lambda_i>0} \lambda_i+\lambda\sum_{\lambda_i<0} \lambda_i
\end{equation*}
and $\lambda_i$ are eigenvalues of $M$.

\begin{re}\label{re.univ}
We say that a constant depending only on $n, \lambda$ and $\Lambda$ is universal.
\end{re}

Next, we introduce the \emph{Pucci's class}.
\begin{defn}\label{d-Sf}
Let $f\in L^{n+1}(\Omega)$. We say that $u\in \underline{S}(\lambda,\Lambda,f)$ (resp. $u\in \overline{S}(\lambda,\Lambda,f)$ ) if $u$ is an $L^{n+1}$-viscosity subsolution (resp., supersolution) of
\begin{equation}\label{SC_Sf}
\begin{aligned}
  &u_t-\mathcal{M}^+(D^2u,\lambda,\Lambda)= f\\
  &\left(\mathrm{resp.}, ~u_t-\mathcal{M}^-(D^2u,\lambda,\Lambda)= f\right).
\end{aligned}
\end{equation}

Define
\begin{equation*}
\begin{aligned}
  S(\lambda,\Lambda,f)= \underline S(\lambda,\Lambda,f)\cap \overline S(\lambda,\Lambda,f)\\
  ~~\mbox{ and }~~
{S^*(\lambda,\Lambda,f)=\underline{S}(\lambda,\Lambda,|f|)\cap \overline{S}(\lambda,\Lambda,-|f|).}
\end{aligned}
\end{equation*}
\end{defn}

\begin{re}\label{re.S}
We remark here that $S(\lambda,\Lambda,f)\subset S^*(\lambda,\Lambda,f)$ and any viscosity solution of \eqref{e.Dirichlet} belongs to $S(\lambda,\Lambda,f+F(0))$ (see \cite[Proposition 3.10]{MR1135923}).
\end{re}

In addition, we introduce the definition of \emph{Dini functions}.
\begin{defn}
Given a one-variable function $\omega:[0,+\infty)\rightarrow [0,+\infty)$, we say that $\omega$ is a Dini function if \\
$(i)$ $\omega$ is nondecreasing;\\
$(ii)$ $\omega$ satisfies the Dini condition
\begin{equation}\label{dini}
	\int_0^{r_0}\frac{\omega(r)}{r}dr<\infty \quad  ~~\mbox{ for some }~~\ r_0>0.
\end{equation}	
\end{defn}
\begin{re}\label{re2.2}
Typical Dini functions are $\omega(r)=r^{\alpha}$ for some $\alpha >0$ and $\omega(r)=(\ln r)^{-\beta}$ for some $\beta>1$.
\end{re}

Now, we introduce some pointwise characterizations of smoothness of functions.

\begin{defn} (Pointwise $C^{0,1}$)
Let $\Omega\subset \mathbb{R}^{n+1}$ be a bounded domain. We say that $f$ is $C^{0,1}$ at $(x_0,t_0)\in \overline{\Omega}$ or $f\in C^{0,1}(x_0,t_0)$ if there exist $K,r_0>0$ such that
\begin{equation}\label{c01}
	|f(x,t)-f(x_0,t_0)|\leq K|(x,t)-(x_0,t_0)|,\quad \forall (x,t)\in \Omega\cap Q_{r_0}(x_0,t_0).
\end{equation}
Then define
$$[f]_{C^{0,1}(x_0,t_0)}:=\min\{K: \eqref{c01}\ holds\ with\ K\}$$
and
$$\|f\|_{C^{0,1}(x_0,t_0)}:=|f(x_0,t_0)|+[f]_{C^{0,1}(x_0,t_0)}.$$	
If   $f\in C^{0,1}(x,t)$ for any $(x,t)\in\Omega$ with the same $r_0$ and
$$\|f\|_{C^{0,1}(\overline\Omega)}:=\|f\|_{L^\infty(\Omega)}+\sup_{(x,t)\in\Omega}[f]_{C^{0,1}(x,t)}<+\infty,$$	
we say that $f\in C^{0,1}(\overline\Omega)$.
\end{defn}

\begin{defn} (Pointwise $C^{1,\mathrm{Dini}}$)\label{de1.1}
Let $\Omega\subset \mathbb{R}^{n+1}$ be a bounded set (not necessarily a domain). We say that $f$ is $C^{1,\mathrm{Dini}}$ at $(x_0,t_0)\in \overline{\Omega}$ or $f\in C^{1,\mathrm{Dini}}(x_0,t_0)$ if there exist $r_0>0$, a vector $l\in \mathbb{R}^n$ and a Dini function $\omega_f$ such that
\begin{equation}\label{c1d}
		|f(x,t)-f(x_0,t_0)-l\cdot (x-x_0)|\leq |(x,t)-(x_0,t_0)|\omega_f(|(x,t)-(x_0,t_0)|),\ \forall (x,t)\in \Omega\cap Q_{r_0}(x_0,t_0).
\end{equation}
Then define
\begin{equation*}
  Df(x_0):=l,\quad [f]_{C^{1,\mathrm{Dini}}(x_0,t_0)}:=I_{\omega}+\omega_f(r_0)
\end{equation*}
and
\begin{equation*}
  \|f\|_{C^{1,\mathrm{Dini}}(x_0,t_0)}:=|f(x_0,t_0)|+|l|+[f]_{C^{1,\mathrm{Dini}}(x_0,t_0)},
\end{equation*}
where $ I_{\omega}$ is the Dini integral defined as
\begin{equation}\label{e1.2}
  I_{\omega}:=\int_{0}^{r_0} \frac{\omega_f (r)}{r} dr.
\end{equation}
If $f\in C^{1,\mathrm{Dini}}(x,t)$ for any $(x,t)\in\Omega$ with the same $r_0$ and
\begin{equation*}
  \|f\|_{C^{1,\mathrm{Dini}}(\overline\Omega)}:=\|f\|_{C^1(\overline\Omega)}
  +\sup_{(x,t)\in\Omega}[f]_{C^{1,\mathrm{Dini}}(x,t)}<+\infty,
\end{equation*}
we say that $f\in C^{1,\mathrm{Dini}}(\overline\Omega)$.
\end{defn}

\begin{re}
If  $\Omega=D\times(0,T]$ with $D\subset \mathbb{R}^n$ being a bounded Lipschitz domain, the above definitions of $C^{0,1}(\overline\Omega)$ and  $C^{1,\mathrm{Dini}}(\overline\Omega)$ are equivalent to the classical definitions (see \cite{MR167862} and \cite{MR1713596}).
\end{re}

\begin{defn} (Pointwise $C^{-1,\mathrm{Dini}}$)\label{de1.2}
	Let $\Omega\subset \mathbb{R}^{n+1}$ be a bounded domain. We say that $f$ is $C^{-1,\mathrm{Dini}}$ at $(x_0,t_0)\in \overline{\Omega}$ or $f\in C^{-1,\mathrm{Dini}}(x_0,t_0)$ if there exist $r_0>0$ and a Dini function $\omega_f$ such that
	\begin{equation}\label{c1}
	\|f\|^*_{L^{n+1}(\Omega\cap Q_r(x_0,t_0))}\leq r^{-1}\omega_f(r),\quad \forall 0<r<r_0,
	\end{equation}
where
\begin{equation*}
	\|f\|^*_{L^{n+1}(\Omega\cap Q_r(x_0,t_0))}:=\left(\frac{1}{|\Omega\cap Q_r(x_0,t_0)|}\displaystyle\int_{\Omega\cap Q_r(x_0,t_0)}|f|^{n+1}\right)^{1/{n+1}}.
\end{equation*}
Then define
$$\|f\|_{C^{-1,\mathrm{Dini}}(x_0,t_0)}:=I_{\omega}+\omega_f(r_0),$$
where $I_{\omega}$ is defined as in \eqref{e1.2}.

If $f\in C^{-1,\mathrm{Dini}}(x,t)$ for any $(x,t)\in\Omega$ with the same $r_0$ and
$$\|f\|_{C^{-1,\mathrm{Dini}}(\overline\Omega)}:=
	\sup_{(x,t)\in\Omega}\|f\|_{C^{-1,\mathrm{Dini}}(x,t)}<+\infty,$$	
	we say that $f\in C^{-1,\mathrm{Dini}}(\overline\Omega)$.
\end{defn}

\begin{re}
By  H${\ddot{o}}$lder's inequality, we deduce that if $f\in L^p(\Omega\cap Q_1)$ for some $p>n+2$, then $f\in C^{-1,\mathrm{Dini}}(0)$.
\end{re}

Next, we introduce the definitions of the geometrical conditions on domains, under which we establish the boundary Lipschitz regularity, the Hopf lemma and the global $W^{2,\delta}$ regularity.

\begin{defn} ($C^{1,\mathrm{Dini}}$ condition on domains)\label{de1.3}
Let $\Omega \subset \mathbb{R}^{n+1}$ be a bounded domain and $(x_0,t_0)\in \partial \Omega$. We say that $\Omega$ satisfies the exterior $C^{1,\mathrm{Dini}}$ condition at $(x_0,t_0)$ if there exist $r_0>0$ and  a new  coordinate system $\{x_1,..,x_n,t\}$ (by rotating and translating with respect to $x$ and only translating with  respect to $t$) such that
\begin{equation}\label{edini}
	(x_0,t_0)=0, \quad
	Q_{r_0}\cap\{(x',x_n,t):x_n< -|(x',t)|\omega_\Omega(|(x',t)|)\}
	\subset Q_{r_0}\cap\Omega^c,
\end{equation}
where $\omega_\Omega$ is a Dini function.

We say that $\Omega$ satisfies the interior $C^{1,\mathrm{Dini}}$ condition at $(x_0,t_0)$ if
\begin{equation}\label{idini}
	(x_0,t_0)=0, \quad
	Q_{r_0}\cap\{(x',x_n,t):x_n> |(x',t)|\omega_\Omega(|(x',t)|)\}
	\subset Q_{r_0}\cap\Omega.
\end{equation}

Let $\Gamma\subset \partial\Omega$ be relatively open. We say that $\Omega$ satisfies the exterior (interior) $C^{1,\mathrm{Dini}}$ condition on $\Gamma$ if $\partial \Omega$ satisfies the exterior (interior) $C^{1,\mathrm{Dini}}$ condition at any point of $\Gamma$ with the same
$r_0$ and $\omega_{\omega}$.

We say that $\Omega$ is an exterior (interior) $C^{1,\mathrm{Dini}}$ domain if $\partial \Omega$ satisfies the exterior (interior) $C^{1,\mathrm{Dini}}$ condition on $\Gamma$ for any relatively open
$\Gamma\subset \partial \Omega$ with $\mathrm{dist}(\partial_b \Omega,\Gamma)>0$.
\end{defn}

\begin{re}\label{re2.3}
In above definition, $\partial_b \Omega$ denotes the bottom boundary of $\Omega$ and $\mathrm{dist}$ denotes the parabolic distance.
\end{re}

\begin{re}\label{re.domain}
In this paper, if we use \Cref{de1.1}, \Cref{de1.2} for a function $f$ or \Cref{de1.3} for a domain $\Omega$, we always denote by $\omega_f$ or $\omega_{\Omega}$ the Dini function respectively. In addition, we always assume that $r_0=1$ in above definitions.
\end{re}

\section{Boundary Lipschitz Regularity and the Hopf Lemma}\label{Lip_Hopf}
In this section, we prove the boundary Lipschitz regularity in \Cref{b_Lip} and the Hopf lemma in \Cref{b_Hopf}. The perturbation idea is used for the proofs, i.e., the curved boundary with $C^{1,\mathrm{Dini}}$ conditions can be treated as the perturbation of a hyperplane.

\subsection{Boundary Lipschitz Regularity}\label{b_Lip}~\\

Based on the known $C^{1,\alpha}$ estimates on flat boundaries, the boundary pointwise Lipschitz regularity for $S^*(\lambda,\Lambda,f)$ on the exterior $C^{1,\mathrm{Dini}}$ domain is then established through a perturbation argument. Set
$$ Q_r^+(x_0,t_0)=Q_r(x_0,t_0)\cap \{(x,t):x_n>x_{0,n}\},$$
and
$$ S_r(x_0,t_0)=\{(x',x_{0,n},t):|x'-x_0'|<r,t_0-r^2<t\leq t_0\},$$
where $x_{0,n}$ is the $n$-th component of $x_0$. We denote $Q_r^+(0)$ and $S_r(0)$ as $Q_r^+$ and $S_r$ respectively for simplicity.

We first recall the \emph{boundary pointwise $C^{1,\alpha}$ regularity on flat boundaries} which is established in \cite{MR688919} (see also \cite{MR1139064, LZ_2022}).

\begin{lemma}\label{l3.1}
Let $u$ be a viscosity solution of
	\begin{equation*}
	\left\{
	\begin{aligned}
        u&\in S(\lambda,\Lambda,0)   &\mathrm{in} \ \   &Q_1^+,\\
		u&=0   &\mathrm{on}  \ \  &S_1.
	\end{aligned}
	\right.
\end{equation*}
Then there exists a universal constant $0<\alpha<1$ such that $u\in C^{1,\alpha}(0)$, i.e., there exists a constant a such that
	\begin{equation*}
	|u(x,t)-ax_n|\leq C_1x_n|(x,t)|^{\alpha}
			\|u\|_{L^{\infty}(Q_1^+)},\quad \forall (x,t)\in  Q_{1/2}^+
	\end{equation*}
	and
	\begin{equation*}
		|a|\leq C_1	\|u\|_{L^{\infty}(Q_1^+)},
	\end{equation*}
	where $C_1$ is universal.
\end{lemma}

\begin{re}
Lee and Yun \cite{MR4656656} proved that for any $0<\alpha<1$, we have the boundary $C^{1,\alpha}$ regularity if the ellipticity constants $\Lambda,\lambda$ are sufficiently close.
\end{re}

For convenience, we present the following scaling version of \Cref{l3.1}.

\begin{cor}\label{c3.1}
Let $0<\alpha<1$ and $C_1$ be as in Lemma \ref{l3.1}. Suppose that $u$ is a viscosity solution of
\begin{equation*}
	\left\{
	\begin{aligned} u&\in S(\lambda,\Lambda,0)   &\mathrm{in} \ \   & Q_r^+(x_0,t_0),\\
		u&=0   &\mathrm{on}  \ \  & S_r(x_0,t_0).
	\end{aligned}
	\right.
\end{equation*}
Then  $u\in C^{1,\alpha}(x_0,t_0)$, i.e., there exists a constant $a_r$ such that
\begin{equation}\label{3.4}
	|u(x,t)-a_r(x_n-x_{0,n})|\leq \frac{C_1}{r^{1+\alpha}}(x_n-x_{0,n})|(x-x_0,t-t_0)|^{\alpha}
	\|u\|_{L^{\infty}( Q_r^+(x_0,t_0))},\ \forall (x,t)\in \  Q_{r/2}^+(x_0,t_0)
\end{equation}
and
\begin{equation*}
	|a_r|\leq \frac{C_1}{r}	\|u\|_{L^{\infty}( Q_r^+(x_0,t_0))}.
\end{equation*}
\end{cor}

\medskip

Next, we give the~\\
\noindent\textbf{Proof of \Cref{t3.1}.} \textbf{Step 1}. Firstly, we make some normalization.
Set
\begin{equation}\label{3.6}
	\omega(r)=\max\{\omega_\Omega(r),\omega_f(r),\omega_g(r)\}.
\end{equation}
Since $\Omega$, $f$ and $g$ satisfy the $C^{1,\mathrm{Dini}}$ or the $C^{-1,\mathrm{Dini}}$ conditions, we may assume that
$$\omega(1)\leq \eta^2,\quad \int_0^{1}\frac{\omega(r)}{r}dr\leq \eta^2,$$
where $0<\eta<1/8$ is universal (to be specified later).
Furthermore, we may assume that
$$u(0)=g(0)=0,\quad Dg(0)=0.$$
Otherwise, we consider $u-g(0)-Dg(0)\cdot x$.
For simplicity, denote by
\begin{equation*}
K=\|u\|_{L^{\infty}(\Omega\cap Q_1)}+\|f\|_{C^{-1,\mathrm{Dini}}(0)}+\|g\|_{C^{1,\mathrm{Dini}}(0)}.
\end{equation*}

\textbf{Step 2}. Secondly, we show the following claim:\\
\textbf{Claim.} There exist positive constants $\bar C$, $\widetilde C$ and a nonnegative sequence $a_k$ ($k\geq -1$) such that for any $k=0,1,2,...$,
we have
\begin{equation}\label{3.7}
	\sup_{\Omega_{\eta^k}}(u-a_k x_n)\leq \widetilde C  K \eta^k A_k
\end{equation}
and
\begin{equation}\label{3.8}
	|a_k-a_{k-1}|\leq \bar C \widetilde C K A_k
\end{equation}
with
\begin{equation}\label{3.9}
A_0=\eta^2,\quad A_{k+1}=\max\{\omega(\eta^k),\eta^{\alpha/2}A_{k}\},
\end{equation}
where $0<\alpha<1$ is as in \Cref{l3.1} and $\bar C$, $\widetilde C$ are universal.

\medskip

 For $k=0$,  by setting $a_{-1}=a_0=0$ and taking
$\widetilde C$ large enough such that
\begin{equation}\label{3.11}
\widetilde C A_0=\widetilde C \eta^2 \geq 1,
\end{equation}
then \eqref{3.7} and \eqref{3.8} hold for  $k=0$.
We suppose that \eqref{3.7} and \eqref{3.8} hold for $k$ and prove that they hold
for $k+1$.

Set
$$r=\eta^k/2,\quad \Omega_r=\Omega\cap Q_r, \quad (\partial\Omega)_r=\partial\Omega\cap Q_r,$$
$$\widetilde Q_r^+:=Q_r^+(0,-r\omega(r),0)=Q_r^+-r\omega(r)(e_n,0), \quad \widetilde \Omega_r=\Omega\cap \widetilde Q_r^+$$
and
$$\widetilde S_r:=S_r(0,-r\omega(r),0)=S_r-r\omega(r)(e_n,0).$$
That is, we consider the domain $\widetilde Q_r^+$ which is a translation of $Q_r^+$ by $r\omega(r)$ along $e_n$. The domain $\widetilde Q_r^+$ contains the boundary $\partial \Omega$ near the origin $0$.
Since $\omega(\eta)\leq\omega(1)\leq \eta^2\leq 1/8$, we have
$$\Omega_{r/2}\subset \widetilde \Omega _r\subset \Omega_{2r}.$$

In the following, we prove \eqref{3.7} and \eqref{3.8} for $k+1$ by the idea of perturbation. Inspired by \cite{MR4713521}, we regard $f$, $g$ and $\partial\Omega$  as perturbation of $0$, $0$ and a hyperplane respectively.

Let $v$ be a viscosity solution of
	\begin{equation}\label{v1}
	\left\{
	\begin{aligned}
		&v_t-\mathcal{M}^+(D^2 v,\lambda,\Lambda)=0   &\mathrm{in} \ \   &\widetilde Q_r^+,\\
		&v=0\   &\mathrm{on}   \ \ &\widetilde S_r,\\
		&v=\widetilde C K \eta^k A_k     &\mathrm{on}  \ \  &\widetilde Q_r^+\backslash \widetilde S_r.
	\end{aligned}
	\right.
\end{equation}
According to the maximum principle, we have
\begin{equation}\label{v-max}
0 \leq v\leq \widetilde C K \eta^k A_k~~\mbox{ in }~~ \widetilde Q_r^+.
\end{equation}
By \Cref{c3.1}, there exists $a_r$ such that for any $(x,t)\in \Omega_{2\eta r} \subset \widetilde{Q}^+_{r/2}$,
\begin{equation}\label{v3}
	\begin{aligned}
		|v-a_r (x_n+r\omega(r))|
		&\leq \frac{C_1}{r^{1+\alpha}}(x_n+r\omega(r))|(x+r\omega(r)e_n,t)|^{\alpha}
		\|v\|_{L^{\infty}(\widetilde Q_r^+)},\\
		&\leq \frac{C_1}{r^{1+\alpha}}|(x+r\omega(r)e_n,t)|^{1+\alpha}
		\widetilde C K \eta^{k} A_{k}
	\end{aligned}
\end{equation}
and
\begin{equation}\label{v4}
	\begin{aligned}
		|a_r|
		\leq \frac{C_1}{r}	\|v\|_{L^{\infty}(\widetilde Q_r^+)}
		\leq \frac{C_1}{r}	\widetilde C K \eta^{k} A_{k}.
	\end{aligned}
\end{equation}
For any $(x,t)\in \Omega_{\eta^{k+1}}=\Omega_{2\eta r}=\Omega\cap Q_{2\eta r}$, recalling $0<\eta<1/8$, we have
$|x+r\omega(r)e_n|<2\eta r+\eta^2 r<r/2$. In addition, since $\Omega$ satisfies the exterior $C^{1,\mathrm{Dini}}$ condition at $0\in\partial\Omega$, we have $x_n>-r\omega(r)e_n$. Hence,
\begin{equation*}
  \Omega_{\eta^{k+1}}=\Omega_{2\eta r}\subset \Omega_{r/4}
  \subset \widetilde Q_{r/2}^+
\end{equation*}
and
\begin{equation*}
    |(x+r\omega(r)e_n,t)|
    \leq \left(\left(2\eta r+\eta^2 r\right)^2+(2\eta r)^2\right)^{\frac{1}{2}}
    \leq 3\eta r, \quad \forall (x,t)\in \Omega_{2\eta r}.
\end{equation*}
It then follows from \eqref{3.9} and \eqref{v3} that
\begin{equation}\label{v3'}
	\begin{aligned}
       \|v-a_r (x_n+r\omega(r))\|_{L^{\infty}(\Omega_{\eta^{k+1}})}
        &=\|v-a_r (x_n+r\omega(r))\|_{L^{\infty}(\Omega_{2\eta r})}
		\leq \|v-a_r (x_n+r\omega(r))\|_{L^{\infty}(\widetilde Q_{r/2}^+)}\\
		&\leq \frac{C_1}{r^{1+\alpha}}(3\eta r)^{1+\alpha}
		\widetilde C K \eta^{k} A_{k}
		\leq 3^{1+\alpha} C_1 \eta^{\frac{\alpha}{2}}
		\widetilde C K \eta^{k+1} \eta^{\frac{\alpha}{2}}A_{k}\\
		&\leq 9 C_1 \eta^{\frac{\alpha}{2}}
		\widetilde C K \eta^{k+1} A_{k+1}.
	\end{aligned}
\end{equation}
Recall that $\omega(r)\leq \eta^2<1/8$. Deduce from \eqref{v4} that
\begin{equation}\label{v4'}
	\begin{aligned}
		|a_r r\omega(r)|
		&\leq \frac{C_1}{r}	\widetilde C K \eta^{k} A_{k} r \eta^2
		= C_1  \eta^{1-\frac{\alpha}{2}}	\widetilde C K \eta^{k+1} \eta^{\frac{\alpha}{2}} A_{k} \\
		&\leq C_1  \eta^{1-\frac{\alpha}{2}}	\widetilde C K \eta^{k+1} A_{k+1}.
	\end{aligned}.
\end{equation}
Consequently,
\begin{equation}\label{v2}
	\begin{aligned}
	\|v-a_r x_n\|_{L^\infty (\Omega_{\eta^{k+1}})}
     \leq \left(9\eta^{\frac{\alpha}{2}}+\eta^{1-\frac{\alpha}{2}}\right) C_1  \widetilde C K \eta^{k+1} A_{k+1}.
	\end{aligned}
\end{equation}

\medskip

Set
$$w=u-a_k x_n-v.$$
By noting \eqref{v1} and \eqref{v-max},
$$v=\widetilde C K \eta^k A_k \ \mbox{ on } \ \widetilde Q_r^+\backslash \widetilde S_r
~~\mbox{ and }~~
0 \leq v\leq \widetilde C K \eta^k A_k\ \ \mathrm{in}\ \widetilde Q_r^+.$$ Then $w$ satisfies
	\begin{equation}\label{w}
	\left\{
	\begin{aligned}
		&w\in \underline S(\lambda,\Lambda,|f|)   &\mathrm{in} \ \   &\Omega\cap \widetilde Q_r^+,\\
		&w\leq g-a_k x_n\   &\mathrm{on}   \ \ &\partial\Omega\cap \widetilde Q_r^+,\\
		&w\leq 0    &\mathrm{on}  \ \  &\Omega\cap \partial\widetilde Q_r^+.
	\end{aligned}
	\right.
\end{equation}
By the Alexsandrov-Bakel'man-Pucci-Krylov-Tso maximum principle \cite[Theorem 3.14]{MR1135923},
\begin{equation*}
  \sup_{\Omega_{\eta^{k+1}}} w
	\leq \sup_{\widetilde \Omega_{r}} w
	\leq |\sup_{\partial \Omega_r} (g-a_k x_n)|+C_2 r^{\frac{n}{n+1}}\|f\|_{L^{n+1}(\widetilde \Omega_{r})},
\end{equation*}
where $C_2$ is universal.
Noticing that
$$a_k x_n\geq -a_k r\omega(r),\quad \forall x\in \partial\Omega\cap \widetilde Q_r^+,$$
we deduce that
$$ |\sup_{\partial\Omega_r} (g-a_k x_n)|
\leq\sup_{\partial\Omega_r} |g|+ a_k r\omega(r).$$
Since $g\in C^{1,\mathrm{Dini}}(0)$,
\begin{equation*}
  \begin{aligned}
    \sup_{\partial\Omega_r} |g|
	&\leq \|g\|_{L^\infty (\partial\Omega\cap \widetilde Q_r^+)}
	\leq \|g\|_{L^\infty (\partial\Omega\cap Q_{2r})}
	\leq K\cdot 2r \omega(2r)\\
	&\leq K \eta^k \omega(\eta^k)
	=\frac{1}{\eta^{1+\frac{\alpha}{2}}} K \eta^{k+1} \eta^{\frac{\alpha}{2}}\omega(\eta^k)
	\leq \frac{1}{\eta^{1+\frac{\alpha}{2}}} K \eta^{k+1} A_{k+1}.
\end{aligned}
\end{equation*}

To estimate $a_k r\omega(r)$, we first prove that
\begin{equation}\label{Ak}
	\sum\limits_{k=0}^{\infty} A_k\leq 4\eta^2.
\end{equation}	
Indeed, by \eqref{3.9} and direct calculations,
$$\sum\limits_{k=0}^{\infty} A_k=\sum\limits_{k=1}^{\infty} A_k+A_0
\leq \sum\limits_{k=0}^{\infty}\omega(\eta^k)+\eta^{\alpha/2}\sum\limits_{k=0}^{\infty} A_k+\eta^2$$
and then
$$\begin{aligned}
	\sum\limits_{k=0}^{\infty} A_k
	&\leq \frac{1}{1-\eta^{\alpha/2}} \left(\sum\limits_{k=0}^{\infty} \omega(\eta^k)+\eta^2\right)
	=\frac{1}{1-\eta^{\alpha/2}} \left(\sum\limits_{k=0}^{\infty} \frac{\omega(\eta^k)(\eta^{k}-\eta^{k+1})}{\eta^{k}-\eta^{k+1}}+\eta^2\right) \\
	&\leq \frac{1}{1-\eta^{\alpha/2}} \left(\frac{1}{1-\eta}\int_0^1\frac{\omega(r)}{r}dr+\eta^2\right)
	\leq \frac{(2-\eta)}{(1-\eta^{\alpha/2})(1-\eta)}\eta^2.
\end{aligned}$$
Then we obtain \eqref{Ak} by taking $\eta$ small enough such that
\begin{equation}\label{3.10}
	\frac{(2-\eta)}{(1-\eta^{\alpha/2})(1-\eta)}\leq 4.
\end{equation}	
Hence,
$$\begin{aligned}
    a_k r\omega(r)
	&\leq \eta^{k} \omega(\eta^k)\sum_{i=0}^{k}|a_i-a_{i-1}|
	\leq \eta^{k} \omega(\eta^k)\bar C\widetilde C K \sum_{i=0}^{k} A_i\\
	&\leq 4\eta^2 \eta^{k} \omega(\eta^k)\bar C\widetilde C K
	= 4 \eta^{1-\frac{\alpha}{2}} \bar C \widetilde C K \eta^{k+1}\eta^{\frac{\alpha}{2}}\omega(\eta^k)
	\leq 4 \eta^{1-\frac{\alpha}{2}} \bar C \widetilde C K  \eta^{k+1}A_{k+1}.
\end{aligned}$$
Since $f\in C^{-1,\mathrm{Dini}}(0)$,
$$\begin{aligned}
\|f\|_{L^{n+1}(\widetilde \Omega_{r})}
	\leq \|f\|_{L^{n+1}(\Omega_{2r})}
	\leq Cr^{\frac{n+2}{n+1}} \|f\|^*_{L^{n+1}(\Omega_{2r})}
	\leq C r^{\frac{n+2}{n+1}}r^{-1}K\omega(2r)
	\leq Cr^{\frac{1}{n+1}}K\omega(2r),
\end{aligned}$$
where $C$ depends only on $n$. It then follows that
\begin{equation*}
  \begin{aligned}
C_2 r^{\frac{n}{n+1}} \|f\|_{L^{n+1}(\widetilde \Omega_{r})}
\leq &C_2r^{\frac{n}{n+1}}C(n)r^{\frac{1}{n+1}}K\omega(2r)
\leq C_3  Mr \omega(2r)\\
\leq &C_3 K \eta^{k}\omega(\eta^k)
=C_3 K \frac{1}{\eta^{1+\frac{\alpha}{2}}} \eta^{k+1} \eta^{\alpha/2}\omega(\eta^k)
\leq \frac{C_3}{\eta^{1+\frac{\alpha}{2}}}K\eta^{k+1} A_{k+1},
  \end{aligned}
\end{equation*}
where $C_3$ is universal.
Coupling the above estimates, we obtain that
\begin{equation}\label{w2}
\begin{aligned}
	\sup_{\Omega_{\eta^{k+1}}} w
	&\leq \|g\|_{L^\infty (\partial\Omega\cap \widetilde Q_r^+)}+a_k r \omega(r)+C_2 r^{\frac{n}{n+1}}\|f\|_{L^{n+1}(\widetilde \Omega_{r})}\\
	&\leq \frac{1}{\eta^{1+\frac{\alpha}{2}}} K \eta^{k+1} A_{k+1}
	+4 \eta^{1-\frac{\alpha}{2}} \bar C \widetilde C K  \eta^{k+1}A_{k+1}
	+\frac{C_3}{\eta^{1+\frac{\alpha}{2}}}K\eta^{k+1} A_{k+1}
	\\
    & \leq \left(\frac{1+C_3}{\widetilde{C}\eta^{1+\frac{\alpha}{2}}}+
    4 \eta^{1-\frac{\alpha}{2}} \bar C\right)
    \widetilde{C}K\eta^{k+1} A_{k+1}.	
\end{aligned}
\end{equation}
By taking $\eta$ small enough, we have
\begin{equation*}
    9 C_1  \eta^{\frac{\alpha}{2}}<\frac{1}{4}, \quad 9 C_1  \eta^{1-\alpha}<\frac{1}{4}
\end{equation*}
and then by taking $\widetilde{C}$ large enough, we get
\begin{equation*}
  \frac{1+C_3}{\widetilde{C}\eta^{1+\frac{\alpha}{2}}}<\frac{1}{4}.
\end{equation*}

Let $$a_{k+1}=a_k+a_r.$$
By \eqref{v4} and choosing
\begin{equation*}
  \bar C=\frac{2 C_1}{\eta^{\alpha/2}},
\end{equation*}
we have
\begin{equation*}
  |a_{k+1}-a_k|=|a_r|
\leq \frac{C_1}{r}	\widetilde C K \eta^{k} A_{k}
\leq \frac{2 C_1}{\eta^{\alpha/2}}\widetilde C K  A_{k+1}
=\bar{C} \widetilde C K A_{k+1}.
\end{equation*}
It then follows from \eqref{v2} and \eqref{w2} that
\begin{equation*}
	\begin{aligned}
        \sup_{\Omega_{\eta^{k+1}}}(u-a_{k+1} x_n)
		=&\sup_{\Omega_{\eta^{k+1}}}(u-a_kx_n-v+v-a_r x_n)
        =\sup_{\Omega_{\eta^{k+1}}}(w+v-a_r x_n)\\
        \leq& \left(\frac{1+C_3}{\widetilde{C}\eta^{1+\frac{\alpha}{2}}}
        +4 \eta^{1-\frac{\alpha}{2}}\bar{C}
        +9C_1 \eta^{\frac{\alpha}{2}}+ C_1  \eta^{1-\frac{\alpha}{2}}\right)\widetilde C K \eta^{k+1} A_{k+1}\\
        \leq &\left(\frac{1+C_3}{\widetilde{C}\eta^{1+\frac{\alpha}{2}}}
        +9 C_1 \eta^{1-\alpha}
        +9C_1 \eta^{\frac{\alpha}{2}}\right)\widetilde C K \eta^{k+1} A_{k+1}\\
		\leq& 	\widetilde C K \eta^{k+1} A_{k+1}.
	\end{aligned}
\end{equation*}
That is, \eqref{3.7} and \eqref{3.8} hold for $k+1$. By induction, we conclude the proof of the claim.

\medskip

\textbf{Step 3}. Finally, we show that the above claim implies $u\in C^{0,1}(0)$.
By \eqref{3.8} and \eqref{Ak},
\begin{equation*}
  \sum\limits_{k=0}^{\infty} |a_k-a_{k-1}|\leq 4\eta^2 \bar C \widetilde C K.
\end{equation*}
Clearly, $a_{k}$ converges to some $a$ and
\begin{equation*}
  |a_k-a|\leq \sum\limits_{i=k}^{\infty} |a_i-a_{i+1}|\leq 4\eta^2 \bar C \widetilde C K.
\end{equation*}
Deduce from \eqref{3.7} that
\begin{equation*}
  \sup_{\Omega_{\eta^k}}(u-a x_n)
\leq \sup_{\Omega_{\eta^k}}(u-a_k x_n)+\sup_{\Omega_{\eta^k}}(a_k-a) x_n
\leq  \widetilde C K \eta^k A_k+  4\eta^2 \bar C \widetilde C K \eta^k\leq C K \eta^k,
\end{equation*}
which implies that
\begin{equation*}
  \sup_{\Omega_r}(u-a x_n)\leq CMr,\ \forall r\leq \eta.
\end{equation*}
Similarly, we have $\inf_{\Omega_r}u\geq -CMr$ and the proof is omitted.
Therefore, we obtain
\begin{equation*}
  \|u\|_{L^\infty(\Omega_r)}\leq CMr.
\end{equation*}
~\qed

\subsection{The Hopf Lemma}\label{b_Hopf}~\\

Similarly, based on the Hopf lemma for $S(\lambda,\Lambda,0)$ on flat boundaries, \Cref{t3.1'} on interior $C^{1,\mathrm{Dini}}$ domains are established through the perturbation argument. We first state the \emph{Hopf lemma on flat boundaries}. We refer to \cite[Lemma 2.1]{LZ_2022} for a proof.

\begin{lemma}\label{l3.1'}
	Let $u\geq 0$ be a viscosity solution of
	\begin{equation*}\label{3.3}
		\begin{cases}
            u\in S(\lambda,\Lambda,0)   \ &\mathrm{in} \ \   Q_1^+,\\
			u=0   &\mathrm{on}   \ \  S_1,
		\end{cases}
	\end{equation*}
with $u(e_n/2,-3/4)\geq 1$. Then
	\begin{equation*}
		u(x,t)\geq c_1x_n,\quad \forall (x,t)\in  Q_{1/2}^+,
	\end{equation*}
where $0<c_1<1/2$ is universal.
\end{lemma}

\medskip

Then we give the~\\
\noindent\textbf{Proof of \Cref{t3.1'}.} \textbf{Step 1}. First, we make some normalization. Let
\begin{equation*}
  \omega(r)=\max(\omega_{\Omega}(r),\hat{C}\omega_f(r)),\quad \forall ~0<r<1,
\end{equation*}
where $\hat{C}$ is a universal constant to be specified later. Without loss of generality, we assume
\begin{equation}\label{e2.2}
u(e_n/2,-3/4)\geq 1, \quad \omega(1)\leq \eta^2,\quad \int_0^{1}\frac{\omega(r)}{r}dr\leq \eta^2,
\end{equation}
where $0<\eta<1/8$ is a small universal constant to be specified later.

Otherwise, we can do the following normalization, which is more complicated than the boundary Lipschitz regularity and hence we give the details. Since $\Omega$ satisfies the interior $C^{1,\mathrm{Dini}}$  condition at $0\in\partial\Omega$, there exists $0<r_1<1/4$ depending on $n,\lambda,\Lambda$ and $\omega_\Omega$  such that
\begin{equation}\label{e2.5}
\omega_{\Omega}(r_1)\leq \eta^2,\quad \int_0^{r_1}\frac{\omega_{\Omega}(r)}{r}dr\leq \eta^2,
\end{equation}
where $\eta$ is as in \eqref{e2.2}. By the Harnack inequality and the assumption on $f$
$$\|f\|_{C^{-1,\mathrm{Dini}}(0)}\leq \varepsilon_0 u(e_n/2,-1/2),$$
there exists $0<c<1$ depending on $n,\lambda,\Lambda$ and $\omega_\Omega$ such that
\begin{equation}\label{nor}
 u\left(\frac{r_1}{2} e_n,-\frac{3r_1}{4}\right)\geq c u\left(\frac{1}{2} e_n,-\frac{1}{2}\right)-\|f\|_{L^{n+1}(\Omega_1)}
 \geq (c-\varepsilon_0) u\left(\frac{1}{2} e_n,-\frac{1}{2}\right)
 \geq \frac{c}{2}u\left(\frac{1}{2} e_n,-\frac{1}{2}\right)
\end{equation}
by taking $\varepsilon_0\leq c/2$.
Write for simplicity
\begin{equation*}
(x_0,t_0)=\left(\frac{r_1}{2} e_n,-\frac{3r_1}{4}\right).
\end{equation*}

Next, we consider the following transformation:
\begin{equation*}
 y=\frac{x}{r_1}, \quad s=\frac{t}{r_1^2}, \quad
 \tilde{u}(y,s)=\frac{u(x,t)}{u(x_0,t_0)}.
\end{equation*}
One can easily see that
\begin{equation*}
  \tilde{u}\in S^*(\lambda,\Lambda,\tilde{f})\quad  \mathrm{in}\ \widetilde \Omega_1,
\end{equation*}
where
\begin{equation*}
\tilde f(y,s)=\frac{r_1^2 f(x,t)}{u(x_0,t_0)}, \quad
\widetilde \Omega=\frac{\Omega}{r_1}.
\end{equation*}
Clearly, we have
\begin{equation}\label{e2.1}
\tilde{u}(e_n/2,-3/4)=1.
\end{equation}
For any $0<r<1$,
\begin{equation*}
  \begin{aligned}
\|\tilde{f}\|^*_{L^{n+1}(\tilde{\Omega}_r)}
=\frac{r_1^2}{u(x_0,t_0)}\|f\|^*_{L^{n+1}(\Omega_{r_1r})}
\leq \frac{r_1^2}{u(x_0,t_0)} (r_1r)^{-1}\omega_f(r_1r)
=\frac{r_1}{u(x_0,t_0)} r^{-1}\omega_f(r_1r):=r^{-1}\omega_{\tilde{f}}(r).
  \end{aligned}
\end{equation*}
Thus, by combining with \eqref{nor} and the assumption
\begin{equation*}
\|f\|_{C^{-1,\mathrm{Dini}}(0)}\leq \varepsilon_0 u(e_n/2,-1/2),
\end{equation*}
we have
\begin{equation}\label{e2.4-0}
\hat{C}\omega_{\tilde{f}}(1)=\frac{\hat{C}r_1}{u(x_0,t_0)}\omega_f(r_1)
\leq \frac{\hat{C}\varepsilon_0u(e_n/2,-1/2)}{u(x_0,t_0)}\leq 2\varepsilon_0 \hat{C} c^{-1} \leq \eta^2
\end{equation}
by taking $\varepsilon_0$ small enough. Similarly,
\begin{equation}\label{e2.4}
\int_{0}^{1} \frac{\hat{C}\omega_{\tilde{f}}(r)}{r}dr=
\frac{\hat{C}r_1}{u(x_0,t_0)}\int_{0}^{r_1} \frac{\omega_{f}(s)}{s}ds
\leq \frac{\hat{C}\varepsilon_0 u(e_n/2,-1/2)}{u(x_0,t_0)}\leq \eta^2.
\end{equation}

Finally, let us consider $\tilde{\Omega}$. Define
\begin{equation*}
\omega_{\tilde{\Omega}}(r)=\omega_{\Omega}(r_1r),\quad \forall ~0<r<1.
\end{equation*}
By the definition of $\tilde{\Omega}$, we have
\begin{equation*}
  \begin{aligned}
Q_1\cap\{(y,s):y_n\geq |(y',s)|\omega_{\tilde \Omega}(|(y',s)|)\}
=&r_1^{-1}\cdot  Q_{r_1}\cap\{(x,t):x_n\geq |(x',t)|\omega_{\Omega}(|(x',t)|)\}\\
\subset& r_1^{-1}\cdot  Q_{r_1}\cap \Omega=Q_1\cap \tilde{\Omega}.
  \end{aligned}
\end{equation*}
Similarly, we have
\begin{equation*}
  \begin{aligned}
Q_1\cap\{(y,s):y_n\leq -|(y',s)|\omega_{\tilde \Omega}(|(y',s)|)\}
\subset Q_1\cap \tilde{\Omega}^c.
  \end{aligned}
\end{equation*}
By noting \eqref{e2.5},
\begin{equation}\label{e2.6}
\omega_{\tilde{\Omega}}(1)=\omega_{\Omega}(r_1)\leq \eta^2, \quad
\int_{0}^{1} \frac{\omega_{\tilde{\Omega}}(r)}{r}dr
=\int_{0}^{r_1}\frac{\omega_{\Omega}(s)}{s}ds\leq \eta^2.
\end{equation}

From \eqref{e2.1}, \eqref{e2.4-0}, \eqref{e2.4} and \eqref{e2.6}, we conclude that the desired assumption \eqref{e2.2} holds for $\tilde{u}$, $\tilde{f}$ and $\tilde{\Omega}$. Therefore, in the following, we always assume \eqref{e2.2} for $u,f$ and $\Omega$.

\medskip

\textbf{Step 2}. Secondly, we set $$\Omega_r^+=\Omega\cap Q_r^+$$ and  show the following claim.\\
\textbf{Claim.} There exist two universal constants $\bar C\geq 1$, $\tilde a$ and a nonnegative sequence $a_k$ $(k\geq -1)$ such that for any $k=0,1,2,...$,
we have
\begin{equation}\label{3.7'}
	\inf_{\Omega^+_{\eta^{k+1}}}(u-\tilde a x_n+a_k x_n)\geq - \eta^k A_k,
\end{equation}
\begin{equation}\label{3.8'}
	a_k-a_{k-1}\leq \bar C  A_k
\end{equation}
and
\begin{equation}\label{3.81'}
	a_k\leq \tilde a/2,
\end{equation}
where
\begin{equation*}
	A_0=\eta^2,\quad A_{k+1}=\max\{\omega(\eta^k),\eta^{\alpha/2}A_{k}\}
\end{equation*}
and $\alpha$ is defined as in \Cref{l3.1}.

\medskip

For $k=0$, by \Cref{l3.1'}, there exist $\eta_1$ and $0<c_1<1/2$ depending only on $n$, $\lambda$ and $\Lambda$ such that
 $$u(x,t)\geq c_1(x_n-\omega(1))\quad \mathrm{in}\ \ Q^+_{\eta_1}+\omega(1)e_n.$$
Note that $u\geq 0$ and hence
 $$u(x,t)\geq c_1(x_n-\omega(1))\quad \mathrm{in}\ \ \Omega^+_{\eta_1}.$$
Then by setting
$$\tilde a=c_1,\quad a_{-1}=a_0=0$$ and
$\eta\leq \eta_1$, we deduce that
$$\inf_{\Omega_\eta^+}(u-c_1 x_n)
\geq \inf_{\Omega_{\eta_1}^+}(u-c_1 x_n)
\geq -c_1 \omega(1)\geq -c_1 \eta^2\geq -\eta^2.$$
Hence \eqref{3.7'}, \eqref{3.8'} and \eqref{3.81'} hold for $k=0$.

Suppose that they hold for $k$ and we will prove them for $k+1$. Let $r=\eta^{k+1}$ and $v$ be a viscosity solution of
\begin{equation*}
	\left\{
	\begin{aligned}
		&v_t-\mathcal{M}^-(D^2 v,\lambda,\Lambda)=0   &\mathrm{in} \ \   &Q_r^+;\\
		&v=0\   &\mathrm{on}   \ \ & S_r;\\
		&v=- \eta^k A_k     &\mathrm{on}  \ \  &\partial Q_r^+\backslash  S_r.
	\end{aligned}
	\right.
\end{equation*}
By \Cref{c3.1}, there exists $a_r$ such that
\begin{equation}\label{vh2}
	0<a_r\leq \frac{C_1}{r}\|v\|_{L^\infty(Q_{r}^+)}
	\leq \frac{C_1}{r}\eta^k A_k
	= \frac{C_1}{\eta} A_k
\end{equation}
and
\begin{equation}\label{vh3}
  \begin{aligned}
\|v+a_r x_n\|_{L^\infty(Q_{\eta^{k+2}}^+)}	
    =&\|v+a_r x_n\|_{L^\infty(Q_{\eta r}^+)}
	\leq \frac{C_1}{r^{1+\alpha}}(\eta r)^{1+\alpha}\|v\|_{L^\infty(Q_{r}^+)}\\
	\leq& C_1\eta^{1+\alpha}\eta^k A_k
	\leq C_1 \eta^{\alpha/2} \eta^{k+1} A_{k+1},
  \end{aligned}
\end{equation}
where $\alpha$ and $C_1$ are as in Lemma \ref{l3.1}.

Let
$w=u-\tilde a x_n+a_k x_n-v$
Then $w$ satisfies
\begin{equation*}
	\left\{
	\begin{aligned}
		&w_t-\mathcal{M}^-(D^2 w,\lambda,\Lambda)\geq f  &\mathrm{in} \ \   &\Omega_r^+;\\
		&w\geq -\tilde a x_n+a_k x_n\   &\mathrm{on}   \ \ &\partial\Omega\cap Q_r^+;\\
		&w\geq 0    &\mathrm{on}  \ \  &\partial  Q_r^+\cap \bar \Omega.
	\end{aligned}
	\right.
\end{equation*}
By the Alexsandrov-Bakel'man-Pucci-Krylov-Tso maximum principle \cite[Theorem 3.14]{MR1135923},
$$\inf_{\Omega^+_{\eta^{k+2}}} w
\geq \inf_{\Omega^+_{r}} w
\geq -\tilde a r\omega(r)-C_2 r^{\frac{n}{n+1}} \|f\|_{L^{n+1}(\Omega^+_{r})},$$
where $C_2$ is universal. Calculate  that
$$\tilde a r\omega(r)= c_1 \eta^{k+1}\omega(\eta^{k+1})
\leq c_1 \eta^{k+1}A_{k+1}
\leq \frac{1}{2}\eta^{k+1}A_{k+1}.$$
By taking $\hat{C}$ large enough,
\begin{equation*}
  \begin{aligned}
C_2 r^{\frac{n}{n+1}} \|f\|_{L^{n+1}(\Omega^+_{r})}
\leq \frac{1}{4}\hat{C}r^2 \|f\|^*_{L^{n+1}(\Omega^+_{r})}
\leq \frac{1}{4}\hat{C}r\omega_f(r)\leq\frac{1}{4} r\omega(r)\leq \frac{1}{4} \eta^{k+1}A_{k+1}.
  \end{aligned}
\end{equation*}
Then,
\begin{equation}\label{wh}
\inf_{\Omega^+_{\eta^{k+2}}} w	
\geq -\frac{3}{4}\eta^{k+1}A_{k+1}
\end{equation}

Set $$a_{k+1}=a_k+a_r.$$
By taking $\eta$ small enough such that
$$\eta\leq \eta_1, \quad C_1\eta^{\alpha/2}\leq 1/4,$$
we derive from \eqref{vh3} and \eqref{wh} that
\begin{equation*}
	\begin{aligned}
		\inf_{\Omega_{\eta^{k+2}}^+}(u-\tilde a x_n+a_{k+1}x_n)
		&=\inf_{\Omega_{\eta^{k+2}}^+}(u-\tilde a x_n+a_{k}x_n-v+v+a_r x_n)\\
		&=\inf_{\Omega_{\eta^{k+2}}^+}(w+v+a_r x_n)
		\geq  -\eta^{k+1} A_{k+1}.
	\end{aligned}
\end{equation*}
and from \eqref{vh2},
$$a_{k+1}-a_k=a_r
\leq \frac{C_1}{\eta} A_{k}
\leq \frac{C_1}{\eta^{1+\frac{\alpha}{2}}} A_{k+1}
:=\bar C A_{k+1}.$$
In addition, recall that $\sum\limits_{i=0}^{\infty} A_i\leq 4\eta^2$ (see \eqref{Ak}).
Then by choosing $\eta$ small enough we have
$$a_{k+1}\leq \sum\limits_{i=0}^{k+1}|a_i-a_{i-1}|
\leq \bar C \sum\limits_{i=0}^{\infty} A_i\leq 4C_1 \eta^{\frac{\alpha}{2}}\leq \tilde a/2.$$
Consequently, \eqref{3.7'}, \eqref{3.8'} and \eqref{3.81'} hold for $k+1$.

\medskip

\textbf{Step 3}. Finally, we show that the claim in Step 2 implies that
$$u(rl,t)\geq cl_n r,\quad \forall 0<r<\eta, -\eta^2<t\leq 0$$
holds for any $l=(l_1,l_2,...,l_n)\in \mathbb{R}^n$
with $|l|=1$ and $l_n> 0$.
In fact, by \eqref{Ak}, we have $$\lim_{k\rightarrow \infty}A_k=0.$$
Then for any $l\in \mathbb{R}^n$ with $|l|=1$ and $l_n>0$, there exists $k_0\geq 1$ such that for any $k\geq k_0$,
$$ \frac{A_k}{l_n\eta^2}\leq \frac{\tilde a}{4}.$$
Set $$r_0=\eta^{k_0+1}.$$
Then for any $0<r<r_0$ and $-r_0^2<t\leq 0$, there exists $k\geq k_0$ such that
$$\eta^{k+2}\leq r\leq \eta^{k+1}\ \ \mathrm{and}\ \ -\eta^{2(k+1)}<t\leq0.$$
For any $(rl,t)\in \Omega_{\eta^{k+1}}^+$, we derive from  \eqref{3.7'}  that
\begin{equation*}
  u(rl,t)-\tilde a l_n r+a_k l_n r
\geq \inf_{\Omega_{\eta^{k+1}}^+} (u-\tilde a x_n+a_k x_n)
\geq - \eta^k A_k.
\end{equation*}
It then follows from \eqref{3.81'}  that
\begin{equation*}
  u(rl,t)\geq \tilde a l_n r-a_k l_n r- \eta^k A_k
\geq \frac{\tilde a l_n r}{2}-\frac{A_k r}{\eta^2}
\geq \frac{\tilde a l_n r}{4}.
\end{equation*}
~\qed

\section{\texorpdfstring{$W^{2,\delta}$}{W2,δ} Regularity}\label{W2d}

In this section, by combining the boundary pointwise Lipschitz regularity (\Cref{t3.1}) and the known interior $W^{2,p}$ estimates, we obtain the global $W^{2,\delta}$ estimates for any $0<\delta<1$ by virtue of the Whitney decomposition.

First, let's recall the Whitney decomposition for parabolic domains (see \cite[P. 342]{MR3385162} and \cite[Chapter VI.1]{MR0290095}). We shall use parabolic cubes whose sides are parallel to the axes, i.e., a cube $\mathcal{Q}$ (with center $(x,t)$ and side length $r$) can be written as
\begin{equation*}
\mathcal{Q}=[x-r,x+r]^n\times [t-r^2,t+r^2]\quad\mbox{for some}~~(x,t)\in \mathbb{R}^{n+1},~r>0.
\end{equation*}
Denote the parabolic diameter by $\diam (\mathcal{Q})$. In addition, we use $\mathcal{Q}^*$ to denote the cube with the same center as $\mathcal{Q}$ and parabolic diameter $\diam (\mathcal{Q}^*)=\frac{9}{8}\diam (\mathcal{Q})$. We denote by $|\cdot|$ and $\dist$ the parabolic norm and parabolic distance respectively.
\begin{lemma}\label{leA.1}
Let $E\subset \mathbb{R}^{n+1}$ be a closed set. Then there exist sequences of $X_j\in E$ and cubes $\mathcal{Q}_j$ whose interiors are mutually disjoint, such that
\begin{equation}\label{eA.1}
  \begin{aligned}
&(i)~E^c=\bigcup_{j\geq 1}\mathcal{Q}_j=\bigcup_{j\geq 1}\mathcal{Q}^*_j;\\
&(ii)~ C^{-1} \mathrm{diam}(\mathcal{Q}_j)\leq \dist (\mathcal{Q}_j,E)\leq C\mathrm{diam}(\mathcal{Q}_j),\quad \forall ~j\geq 1;\\
&(iii)~\dist (\mathcal{Q}_j, E)=\dist (\mathcal{Q}_j,X_j),\quad \forall ~j\geq 1;\\
&(iv)~C^{-1} |X-X_j|\leq \mathrm{diam}(\mathcal{Q}^*_j)\leq C|X-X_j|,\quad \forall ~j\geq 1,~\forall ~X\in \mathcal{Q}^*_j;\\
&(v)~|X_0-X_j|\leq C|X_0-X|, \quad \forall ~j\geq 1,~\forall ~X_0\in E, ~\forall ~X\in \mathcal{Q}^*_j;\\
&(vi)~X \mbox{ is contained in at most } N \mbox{ cubes } \mathcal{Q}^*_j, \quad \forall ~X\in E^c,
  \end{aligned}
\end{equation}
where all constants $C$ and $N$ in above inequalities depend only on $n$.
\end{lemma}

Whitney decomposition is an effective tool for obtaining regularity in Sobolev spaces. Cao, Li and Wang \cite{MR2754287} utilized it to prove the optimal weighted $W^{2,p}$ estimates for elliptic equations with non-compatible conditions. Li, the second author and Zhang \cite{MR4631039} utilized it to prove $W^{2,p}$ estimates on $C^{1,\alpha}$ domains. Moreover, the Whitney extension can be obtained based the Whitney decomposition and the characterization of nodal sets need Whitney extension theorems (see \cite{lian2024nodalsets}).

Next, we show the summation of $\diam(\mathcal{Q}_j)$ with Lipschitz boundary is bounded as in \cite{MR4631039}. In the following, we always assume that $0\in \partial \Omega$ and there exists $\varphi \in C^{0,1} (S_1)$ with $\|\varphi\|_{C^{0,1}(S_1)} \leq \mathcal{K}$ such that
\begin{equation*}
  \Omega_1=\{x_n > \varphi (x',t)\}\cap Q_1~~\mbox{ and }~~(\partial \Omega)_1=\{x_n=\varphi(x',t)\}\cap Q_1
\end{equation*}
for some positive constant $\mathcal{K}$. We consider $|\cdot|$ as the parabolic measure, i.e.,
\begin{equation*}
  P_{\delta}(\Omega):=\inf \left\{\sum_{i=1}^{N} s^{n+2}:
  \Omega\subset \bigcup_{i=1}^{N} Q^{i}_s, s\leq \delta \right\},
\end{equation*}
where $Q_s^i=Q_s(x_i,t_i)$ and
\begin{equation*}
|\Omega|:=\lim_{\delta \to 0} P_{\delta}(\Omega).
\end{equation*}

\begin{lemma}\label{le.4.2}
  For any $(x_0,t_0)\in \overline{\Omega}_1$ and $r>0$ with $\Omega_r(x_0,t_0)\subset \Omega_1$,
  \begin{equation}\label{e.Mea}
    |\Omega_r(x_0,t_0)\cap \{\dist ((x,t),(\partial \Omega)_1)\leq d\}|\leq C r^{n+1}d \quad~~\mbox{ for }~~d>0,
  \end{equation}
  where $C$ depends only on $n$ and $\mathcal{K}$.
\end{lemma}

\begin{proof}
By recalling $\Omega_r(x_0,t_0)=\Omega \cap Q_r(x_0,t_0)$, we observe that
\begin{equation*}
  \begin{aligned}
  &\Omega_r(x_0,t_0)\cap \{\dist ((x,t),(\partial \Omega)_1)\leq d\}\\
  \subset & \{|x'-x'_0|\leq r, |t-t_0|\leq r^2, \varphi(x',t)\leq x_n \leq \varphi (x',t)+(\mathcal{K}+1)d\}\\
  \leq & Cr^{n+1} d,
  \end{aligned}
\end{equation*}
where $C$ depends only on $n$ and $\mathcal{K}$. Hence, we have \eqref{e.Mea}.
\end{proof}

Denote $d_j:=\mathrm{diam}(\mathcal{Q}_j)$ and set $$\mathcal{F}^k=\bigcup_{j} \{\mathcal{Q}_j:2^{-k-1}<d_j\leq 2^{-k}\}, \quad k=2,3,\dots.$$

\begin{lemma}\label{l.bound.q}
  If $q>n+1$, then
  \begin{equation*}
    \sum_{\widetilde{\mathcal{Q}}_j \subset \Omega_{1/4}} d_j^q \leq C,
  \end{equation*}
  where $C$ depends only on $n, q$ and $\mathcal{K}$.
\end{lemma}

\begin{proof}
For any $\mathcal{Q}_j$, by (ii) and (iii) of \Cref{leA.1} there exists $(y_j,s_j)\in (\partial \Omega)_1$ such that
\begin{equation*}
  \dist (\mathcal{Q}_j,(y_j,s_j))=\dist (\mathcal{Q}_j,(\partial \Omega)_1)\leq C d_j \leq C 2^{-s}.
\end{equation*}
Then
\begin{equation*}
  \dist ((x,t),(\partial \Omega)_1)\leq d_j+\dist (\mathcal{Q}_j,y_j) \leq (C+1) 2^{-k},\quad
  \forall (x,t)\in \mathcal{Q}_j ~~\mbox{ and }~~\mathcal{Q}_j \in \mathcal{F}^k.
\end{equation*}
It follows that
\begin{equation*}
  \mathcal{F}^k\subset \Omega_{1/4} \cap \{\dist ((x,t),(\partial \Omega)_1)\leq (C+1)2^{-k}\}.
\end{equation*}
By \Cref{le.4.2}, we have
\begin{equation*}
  |\mathcal{F}^k|\leq C2^{-k},
\end{equation*}
where $C$ depends only on $n$ and $K$.
Note that
\begin{equation*}
  \bigcup_{\widetilde{\mathcal{Q}}_j\subset \Omega_{1/4}} \mathcal{Q}_j
  =\bigcup_{k=2}^{\infty} \bigcup_{\mathcal{Q}_j\subset \mathcal{F}^k} \mathcal{Q}_j.
\end{equation*}
If $q>n+1$, we deduce that
\begin{equation*}
  \begin{aligned}
  \sum_{\widetilde{\mathcal{Q}}_j}  d_j^q
  \leq & \sum_{k=2}^{\infty} \left(\sum_{\mathcal{Q}_j \in \mathcal{F}^k} (d_j^{q-n-2}\cdot d_j^{n+2})\right)
  \leq \sum_{k=2}^{\infty} \left(2^{-k(q-n-2)}\cdot  (\sum_{\mathcal{Q}_j \in \mathcal{F}^k} d_j^{n+2})\right)\\
  \leq & C\sum_{k=2}^{\infty} 2^{-k(q-n-2)} |\mathcal{F}^k|
  \leq C \sum_{k=2}^{\infty} 2^{-k(q-n-1)} \leq C,
  \end{aligned}
\end{equation*}
where $C$ depends only $n, q$ and $\mathcal{K}$.
\end{proof}

In addition, we introduce the interior \emph{$W^{2,p}$ regularity} proved by Wang \cite[Theorem 5.6]{MR1135923}.
\begin{thm}\label{0}
Let $u$ be a viscosity solution of \eqref{e.Dirichlet}. Suppose that the operator $F$ is convex or concave in $M$ and $f\in L^{p}(Q_1)$ ($n+1\leq p<+\infty$). Then
$u\in W^{2,p}(Q_{1/2})$ and
\begin{equation}\label{1.5}
    \|u\|_{W^{2,p}(Q_{1/2})}\leq C\left(\|u\|_{L^{\infty}(Q_1)}+ \|f\|_{L^{p}(Q_1)}\right),
\end{equation}
where $C$ is universal.
\end{thm}

\begin{re}\label{re.w2p0}
Here we need $p\geq n+1$ because the maximum principle needs at least $f\in L^{n+1}(Q_1)$ (see \cite[Theorem 3.14]{MR1135923}).
\end{re}

\medskip

Now, we show the~\\
\noindent \textbf{Proof of \Cref{t4.1}}.
Let
\begin{equation}\label{h}
K=\|u\|_{L^{\infty}(\Omega_1)}+\|f\|_{L^{p}(\Omega_1)}+ \|g\|_{C^{1,\mathrm{Dini}}((\partial \Omega)_1)}.
\end{equation}
We only need to prove
\begin{equation}\label{l}
\begin{aligned}
\|D^2u\|_{ L^{\delta}(\Omega_{1/12})}+
\|u_t\|_{ L^{\delta}(\Omega_{1/12})}\leq CK,
\end{aligned}
\end{equation}
since \eqref{1.2} can be obtained by a standard covering argument.

First, we claim that
\begin{equation}\label{e.o13}
  \Omega_{r/3}\subset \bigcup_{\mathcal{Q}^*_j \subset \Omega_{r}} \mathcal{Q}_j
\end{equation}
for any $0<r\leq 1$. If not, there  exist  a point $(x,t)\in \Omega_{r/3}$ and a cube $\mathcal{Q}_j$ such that $(x,t)\in \Omega_{r/3}$ but $\mathcal{Q}^*_j \not\subset \Omega_r$. It means that there exists a point $(y,s)\in \mathcal{Q}^*_j$ with $|(y,s)|\geq r$. Then, by (ii) and (iv) of \Cref{leA.1}, we have
\begin{equation*}
  \dist (\mathcal{Q}_j,(\partial \Omega)_1) \geq C d_j
  =\frac{8}{9} C \diam (\mathcal{Q}^*_j)
  \geq \frac{8}{9}(|(y,s)|-|(x,t)|)\geq r/2,
\end{equation*}
which contradicts with
\begin{equation*}
  \dist ((\mathcal{Q}_j,(\partial \Omega)_1)\leq |(x,t)|\leq r/3.
\end{equation*}

Next, we obtain the estimate on every $\mathcal{Q}^*_j$ with the aid of the boundary pointwise Lipschitz regularity. For any $\mathcal{Q}^*_j$, by (ii) and (iii) of \Cref{leA.1} there exists $(y_j,s_j)\in (\partial \Omega)_1$ such that
\begin{equation*}
  \dist (\mathcal{Q}^*_j,(y_j,s_j))=\dist (\mathcal{Q}^*_j,(\partial \Omega)_1)
  \leq \dist (\mathcal{Q}_j,(\partial \Omega)_1)
  \leq C d_j,
\end{equation*}
where $C$ depends only on $n$ and $\omega_{\Omega}$. Then for any $(x,t)\in \mathcal{Q}^*_j$, we have
\begin{equation*}
  |(x,t)-(y_j,s_j)|\leq 2\diam(\mathcal{Q}^*_j)+\dist (\mathcal{Q}^*_j,(y_j,s_j))
  \leq C d_j,
\end{equation*}
where $C$ depends only on $n$ and $\omega_{\Omega}$.

By the boundary $C^{0,1}$ estimates,
\begin{equation}\label{1a}
		\|u-u(y_j,s_j)\|_{L^\infty(\mathcal{Q}^*_j)}\leq  CK d_j,
\end{equation}
where $C$ depends only on $n,\lambda,\Lambda,p,\omega_g$ and $\omega_{\Omega}$.
	
Note that $u-u(y_j,s_j)$ satisfies the same equation as $u$. Then, by H\"{o}lder inequality, \Cref{0} and \eqref{1a}, for any $0<\delta\leq p$,
\begin{equation}\label{impro_d}
		\begin{aligned}
			&\|D^2 u\|_{L^{\delta}(\mathcal{Q}_j)}+\|u_t\|_{L^{\delta}(\mathcal{Q}_j)}\\
          =&\|D^2 (u-u(y_j,s_j))\|_{L^{\delta}(\mathcal{Q}_j)}+\|(u-u(y_j,s_j))_t\|_{L^{\delta}(\mathcal{Q}_j)}\\
           \leq & C d_j ^{(n+2)\left(\frac{1}{\delta}-\frac{1}{p}\right)}
           \left(\|D^2 (u-u(y_j,s_j))\|_{L^{p}(\mathcal{Q}_j)}+\|(u-u(y_j,s_j))_t\|_{L^{p}(\mathcal{Q}_j)}\right)\\
		   \leq & C\left(
			d_j^{\frac{n+2}{\delta}-2}\|u-u(y_j,s_j)\|_{L^\infty(\mathcal{Q}^*_j)}
			+d_j^{(n+2)\left(\frac{1}{\delta}-\frac{1}{p}\right)}
			\|f\|_{L^p(\mathcal{Q}^*_j)}\right)\\
           \leq & CK\left(d_j^{\frac{n+2}{\delta}-1} +d_j^{(n+2)\left(\frac{1}{\delta}-\frac{1}{p}\right)}\right),
		\end{aligned}
	\end{equation}
where $C$ depends on $n,\lambda,\Lambda,p,\omega_g$ and  $\omega_{\Omega}$.
	
Recalling \eqref{e.o13}, we have
\begin{equation*}
  \Omega_{1/{12}}\subset\bigcup_{\mathcal{Q}^*_j\subset \Omega_{1/{4}}} \mathcal{Q}_j.
\end{equation*}
Then sum over $j$,
\begin{equation*}
	\begin{aligned}
        \|D^2u\|_{ L^{\delta}(\Omega_{1/12})}+\|u_t\|_{ L^{\delta}(\Omega_{1/12})}
		\leq & \sum\limits_{\mathcal{Q}^*_j\subset \Omega_{1/4}}
		\left(\|D^2 u\|_{L^{\delta}(\mathcal{Q}_j)}+\|u_t\|_{L^{\delta}(\mathcal{Q}_j)}\right)\\
		\leq &  CK \sum\limits_{\mathcal{Q}^*_k\subset \Omega_{1/4}}
        \left(d_j^{\frac{n+2}{\delta}-1} +d_j^{(n+2)\left(\frac{1}{\delta}-\frac{1}{p}\right)}\right).
	\end{aligned}
\end{equation*}
Since $0<\delta<1$ and $p>n+2$, by \Cref{l.bound.q}, we have \eqref{l}.~\qed

\begin{re}\label{re.w2d1}
By \eqref{impro_d}, we can see that the boundary Lipschitz regularity is used to improve the value of $\delta$ (see also \cite{LI2025109459}).
\end{re}

\bigskip

\printbibliography

\end{document}